\documentclass[aos,MSNbibl,seceqn,dvips]{arximspdf}

\doi{10.1214/13-AOS1115} 
\volume{41}
\issue{3}
\pubyear{2013}
\firstpage{1462}
\lastpage{1484}

\makeatletter

\newproclaim{assump}{Assumption}

\newcommand{\R}{\mathbb{R}}
\newcommand{\E}{\mathbb{E}}
\newcommand{\PP}{\mathbb{P}}

\newcommand{\al}{\alpha}
\newcommand{\la}{\lambda}
\newcommand{\vpi}{\varpi}
\newcommand{\si}{\sigma}
\newcommand{\te}{\theta}
\newcommand{\be}{\beta}
\newcommand{\ep}{\varepsilon}
\newcommand{\De}{\Delta}
\newcommand{\Om}{\Omega}
\newcommand{\ze}{\zeta}

\newcommand{\f}{\mathcal{F}}
\newcommand{\laa}{\mathcal{L}}
\newcommand{\m}{\mathcal{M}}
\newcommand{\n}{\mathcal{N}}
\newcommand{\s}{\mathcal{S}}

\newcommand{\WE}{\widetilde{\mathbb{E}}}
\newcommand{\WV}{\widetilde{V}}
\newcommand{\Wc}{\widetilde{c}}
\newcommand{\Wsi}{\widetilde{\sigma}}
\newcommand{\Wde}{\widetilde{\delta}}
\newcommand{\Wb}{\widetilde{b}}

\newcommand{\Bb}{\overline{b}{}}
\newcommand{\BV}{\overline{V}{}}
\newcommand{\BW}{\overline{W}{}}
\newcommand{\Bh}{\overline{h}{}}
\newcommand{\Bsi}{\overline{\sigma}{}}

\newcommand{\rdn}{\sqrt{\De_n}}

\newcommand{\toop}{\stackrel{\PP}{\longrightarrow}}
\newcommand{\tolfs}{\stackrel{\laa_f-s}{\longrightarrow}}
\newcommand{\tolls}{\stackrel{\laa-s}{\Longrightarrow}}
\newcommand{\tols}{\stackrel{\laa-s}{\longrightarrow}}

\newcommand{\toucp}{\stackrel{\mathrm{u.c.p.}}{\Longrightarrow}}

\newtheorem{cor}{Corollary}[section]
\newtheorem{lem}[cor]{Lemma}
\newtheorem{theo}[cor]{Theorem}
\newproclaim{rema}[cor]{Remark}
\newproclaim{ex}[cor]{Example}

\makeatother

\begin{document}
\begin{frontmatter}

\title{Quarticity and other functionals of volatility: Efficient estimation}
\runtitle{Integrated function of volatility}

\begin{aug}
\author[A]{\fnms{Jean} \snm{Jacod}\ead[label=e1]{jean.jacod@upmc.fr}}
\and
\author[B]{\fnms{Mathieu} \snm{Rosenbaum}\corref{}\ead[label=e2]{mathieu.rosenbaum@upmc.fr}}
\runauthor{J. Jacod and M. Rosenbaum}
\affiliation{Universit\'e Pierre et Marie Curie (Paris 6)}
\address[A]{Institut de Math\'ematiques de Jussieu\\
CNRS---UMR 7586\\
Universit\'e Pierre et Marie Curie\\
4 Place Jussieu\\
75 005 Paris\\
France\\
\printead{e1}} 
\address[B]{Laboratoire de Probabilit\'es\\
\quad et Mod\`eles Al\'eatoires\\
CNRS---UMR 7599\\
Universit\'e Pierre et Marie Curie\\
4 Place Jussieu\\
75 005 Paris\\
France\\
\printead{e2}}
\end{aug}

\received{\smonth{7} \syear{2012}}
\revised{\smonth{3} \syear{2013}}

%
\begin{abstract}
We consider a multidimensional It\^o semimartingale regularly sampled
on $[0,t]$ at high frequency $1/\Delta_n$, with $\Delta_n$ going to
zero. The goal of this paper is to provide an estimator for the
integral over $[0,t]$ of a given function of the volatility matrix. To
approximate the integral, we simply use a Riemann sum based on local
estimators of the pointwise volatility. We show that although the
accuracy of the pointwise estimation is at most $\Delta_n^{1/4}$, this
procedure reaches the parametric rate $\Delta_n^{1/2}$, as it is
usually the case in integrated functionals estimation. After a suitable
bias correction, we obtain an unbiased central limit theorem for our
estimator and show that it is asymptotically efficient within some
classes of sub models.
\end{abstract}

%

\end{frontmatter}

\section{Introduction}\label{sec-I}

Let $X$ be a semimartingale, which is observed at discrete times
$i\De_n$ for $i=0,1,\ldots\,$, over a finite time interval $[0,T]$,
with a
discretization mesh $\De_n$ which is small and eventually goes to $0$
(high-frequency setting). One of the main problems encountered in practice
is the estimation of the integrated (squared) volatility (in finance
terms), or equivalently of the continuous part of the quadratic variation
$[X,X]_t$.

By now, this is a well-understood problem, at least when $X$ is an
It\^o semimartingale. For example, in the
continuous one-dimensional case, if $X$ takes the form
\[
X_t=X_0+\int_0^tb_s
\,ds+\int_0^t\si_s
\,dW_s
\]
the approximate quadratic variation
$\sum_{i=1}^{[t/\De_n]}(X_{i\De_n}-X_{(i-1)\De_n})^2$, which of
course converges to
$[X,X]_t=\int_0^t\si^2_s\,ds$, enjoys a central limit theorem (CLT): the
difference between these two processes, normalized by $\frac{1}{\rdn
}$, converges
stably in law to a limit which is conditionally on $X$ a continuous Gaussian
martingale with quadratic variation (equivalently, with variance)
twice the so-called ``quarticity,'' that is, $2\int_0^t\si_s^4\,ds$.

Although later we consider a much more general framework, allowing $X$
to be multi-dimensional and with jumps, in the \hyperref[sec-I]{Introduction} we pursue the
discussion in this special one-dimensional continuous case. In various
statistical problems one needs to estimate not only the quarticity, but
functionals of the form
\[
V(g)_t=\int_0^tg(c_s)
\,ds\qquad \mbox{where } c_s=\si^2_s
\]
(for relatively general test functions $g$, and to derive associated
CLTs, see~\cite{JP}); notice that we plug in the
``spot'' squared volatility $c_t$ rather than $\si_t$, since in any
case it
is impossible to determine the sign of $\si_t$ on the basis of the
observation of the path $t\mapsto X_t$. The case $g(x)=x$ corresponds
to the usual integrated volatility, and $g(x)=x^2$ to the quarticity.

Toward this aim, two methods are currently at hand:
\begin{longlist}[(2)]
\item[(1)] The first one is available if $g(x)=\E(f(U(x)_1,\ldots,
U(x)_k))$ for all \mbox{$x\geq0$}, where the $U(x)_j$'s are independent $\n(0,x)$
variables and $f$ is a continuous function
on $\R^k$, of polynomial growth. Then we know that
%
\begin{eqnarray}
\label{I-1}
U^n(f)_t=\De_n \sum
_{i=1}^{[t/\De_n]-k+1} f \biggl(\frac{\De^n_iX}{\rdn},\ldots,\frac{\De^n_{i+k-1}X}{\rdn
} \biggr)\nonumber\\[-8pt]\\[-8pt]
&&\eqntext{\displaystyle \mbox{where } \De^n_iX=X_{i\De_n}-X_{(i-1)\De_n},}
\end{eqnarray}
converges to $V(g)_t$ in probability, and if $f$ is $C^1$ the rate of
convergence is $1/\rdn$, and in the associated CLT the
limiting conditional variance is $\int_0^tF(c_s)\,ds$ for a suitable
function $F$.

\item[(2)] The second one consists in using estimators for the spot volatility
and approximating the integral $V(g)_t$ by Riemann sums, in which the spot
volatility is replaced by its estimator; that is, we set
%
\begin{equation}
\label{I-2}\quad V^n(g)_t=\De_n\sum
_{i=1}^{[t/\De_n]-k_n+1}g\bigl(\widehat{c}{}^n_i
\bigr)\qquad \mbox{where }\widehat{c}{}^n_i=\frac1{k_n
\De_n}\sum_{j=0}^{k_n-1}\bigl(
\De^n_{i+j}X\bigr)^2
\end{equation}
for an arbitrary sequence of integers such that $k_n\to\infty$ and
$k_n\De_n\to0$. Then one knows that $V^n(g)_t\toop V(g)_t$ (when $g$
is continuous and of polynomial growth). But so far nothing is known
about the rate of convergence of these estimators when $k_n$ goes to infinity
(the situation $k_n=k$ not depending on $n$
is studied in~\cite{MZ} where the rate $1/\sqrt{\De_n}$ is obtained
for power
functions).
\end{longlist}

The first method is quite powerful and gives optimal rates, but the
special form of $g$ puts strong constraints
on this function [e.g., it is $C^\infty$ on $(0,\infty)$, and much more].
To tell the truth, in the one-dimensional case, by far the most useful test
functions $g$ are the
powers $g_p(x)=x^p$ (recall that $x\geq0$ here) for $p>0$, which are
associated as above with $f_p(x)=|x|^{2p}/m_{2p}$, where $m_q$ is the $q$th
absolute moment of $\n(0,1)$. Nevertheless, some functions $g$ of interest
might not be, or not in an obvious way, of this form or, more
generally, linear
combinations of functions of this form. In the multivariate case,
however, with $X$ being $d$-dimensional and thus $U$ above
as well, one typically finds asymptotic variances which are complicated
functions of the $d\times d$-dimensional spot volatility. This is, for instance,
the case when studying multipower variations
for integrated volatility estimation in the presence of jumps; see, for
example,
\cite{JP}. In this situation and more generally for an arbitrary (smooth)
function $g$ on the set $\m^+_d$ of all
$d\times d$ symmetric nonnegative matrices, it is rather a difficult
task in
practice to find an integer $k\geq1$ and a function $f$ on $(\R^d)^k$
such that, for all $x\in\m^+_d$, we have
$g(x)=\E(f(U(x)_1,\ldots,U(x)_k))$, where again the $U(x)_j$'s are
($d$-dimensional) i.i.d. $\n(0,x)$.

In addition, this first method does not provide efficient estimation in
general. To see that, consider the toy example $X_t=\si W_t$, where
$\si$ is
a constant, $c=\si^2$, $\De_n=\frac1n$ and $T=1$. We thus observe the
increments $\De^n_iX$ for $i=1,\ldots,n$, or equivalently the $n$ variables
$Y_i=\De^n_iX/\rdn$. These variables are i.i.d. $\n(0,c)$, so the
asymptotically best
estimators for $c$ (efficient in all possible senses, and also the MLE) are
$\widehat{c}_n=\frac1n \sum_{i=1}^n(Y_i)^2=\sum_{i=1}^n(\De
^n_iX)^2$, with
convergence rate $\sqrt{n}$ and asymptotic variance $2c^2$.
If instead one wants to estimate $c^p$ for some $p\neq1$ in $(0,\infty)$,
one can use $\widehat{c}(p)_n=\frac1{n m_{2p}} \sum_{i=1}^n|Y_i|^{2p}
=\frac{n^{p-1}}{m_{2p}}\sum_{i=1}^n|\De^n_iX|^{2p}$, and the
ordinary central
limit theorem tells us that the rate of convergence is again $\sqrt{n}$,
and the asymptotic variance is $\frac{m_{4p}-m_{2p}^2}{m_{2p}^2} c^{2p}$:
this is exactly what the first method above does. But this is not optimal,
the asymptotically optimal estimators being $(\widehat{c}_n)^p$ (the
MLE again), with
rate $\sqrt{n}$ and asymptotic variance $2p^2c^{2p}$, smaller than
the previous one when $p\neq1$. Now, taking
$(\widehat{c}_n)^p$ is exactly what the second method (\ref{I-2}) does.

The aim of this paper is to develop the second method, and in particular
to provide a central limit theorem, with the rate $1/\rdn$ (as it is usually
the case in a nonparametric setting for integrated functionals estimation;
see, e.g.,~\cite{BR,BM}), and with an
asymptotic variance always smaller than if one uses the first method. This
can be viewed as an extension, in several directions, of the ``block
method'' of Mykland and Zhang in~\cite{MZ}. About
efficiency, and despite the title of the paper, we do not really examine
the question in the general nonparametric or semi-parametric setting
assumed below, since even for the simpler problem of estimating the integrated
volatility, the concept of efficiency is not well
established so far. Instead, we will term as ``efficient'' a procedure
which is efficient in the usual sense for the sub-model consisting in
the toy model $X_t=\si W_t$ above, and efficient in the sense of the Hajek
convolution theorem, for the Markov-type model recently studied by
Cl\'ement, Delattre and Gloter in~\cite{CDG} and of the form
%
\begin{equation}
\label{I-3} dX_t=a(X_t)\,dt+f(t,X_t,Y_t)\,dW_t,\qquad
dY_t=\Bb_t\,dt+\Bsi_td\BW_t,
\end{equation}
where $a,f$ are unknown smooth enough functions and $\Bb,\Bsi$ arbitrary
processes and where the two Brownian motions $W,\BW$ are independent.

This will be done in the multivariate setting and when $X$ possibly
has jumps (upon suitably truncating the increments in (\ref{I-2}) if it
is the case, in the spirit of \mbox{\cite{M1,M2}}), and under the additional
assumptions that $c_t$ itself is an It\^o semimartingale and that,
when $X$ jumps, these jumps are summable, which are exactly the same
assumptions under which the truncated versions of $U^n(f)$ in
(\ref{I-1}) converge with rate $1/\rdn$.

The paper is organized as follows: Section~\ref{sec-SET} is devoted to
presenting the assumptions. Results are given in Section~\ref{sec-RES},
and all proofs are gathered in Section~\ref{sec-P}.

\section{Setting and assumptions}\label{sec-SET}

The underlying process $X$ is $d$-dimensional, and observed at the
times $i\De_n$ for $i=0,1,\ldots\,$, within a fixed interval of interest
$[0,t]$. For any process $Y$ we use the notation $\De^n_iY$ defined in
(\ref{I-1})
for the increment over the $i$th observation interval. We assume that the
sequence $\De_n$ goes to $0$. The precise assumptions on $X$ are as follows:

First, $X$ is an It\^o semimartingale on a filtered space
$(\Omega,\f,(\f_t)_{t\geq0},\PP)$. It can be written in its Grigelionis
form, using a $d$-dimensional Brownian motion $W$ and a Poisson random
measure $\mu$ on $\R_+\times E$, where $E$ is an auxiliary Polish space
and with the (nonrandom) intensity measure
$\nu(dt,dz)=dt\otimes\la(dz)$ for some $\si$-finite measure $\la$ on
$E$,
%
\begin{eqnarray}
\label{S-1} X_t&=&X_0+\int_0^tb_s
\,ds+\int_0^t\si_s
\,dW_{s}\nonumber\\
&&{}+\int_0^t\int
_E \delta(s,z) 1_{\{\|\delta(s,z)\|\leq1\}} (\mu-\nu) (ds,dz)
\\
&&{} +\int_0^t\int_E
\delta(s,z) 1_{\{\|\delta(s,z)\|>1\}} \mu(ds,dz).
\nonumber
\end{eqnarray}
This is a vector-type notation: the process $b_t$ is $\R^d$-valued optional,
the process $\si_t$ is $\R^d\otimes\R^d$-valued optional,
$\delta=\delta(\omega,t,z)$ is a predictable $\R^d$-valued function on
$\Om\times\R_+\times E$ and $\|\cdot\|$ denotes the Euclidean norm on any
finite-dimensional linear space. Besides the measurability requirements above,
and for any $r\in[0,2]$, we introduce the assumption:

\renewcommand{\theassump}{(H-$r$)}
\begin{assump}\label{assumHr}
There are a sequence $(J_n)$ of nonnegative
bounded $\la$-integrable functions on $E$ and a sequence
$(\tau_{n})$ of stopping times increasing to $\infty$, such that
%
\begin{eqnarray}
\label{S-3}
t&<&\tau_n(\omega)\quad\Rightarrow\quad\bigl\|b_t(\omega)\bigr\|
\leq n,\qquad\bigl\|\si_t(\omega)\bigr\|\leq n,\nonumber\\[-8pt]\\[-8pt]
t&\leq&\tau_n(\omega)
\quad\Rightarrow\quad\bigl\|\delta(\omega,t,z)\bigr\|^r\wedge1\leq
J_n(z).\nonumber
\end{eqnarray}
\end{assump}

The spot volatility process $c_t=\si_t\si_t^*$ ($^*$ denotes transpose)
takes its values in the set $\m^+_d$ of all nonnegative symmetric
$d\times d$ matrices. We will indeed suppose that $c_t$ is again an It\^o
semimartingale, and we consider the following assumption:

\renewcommand{\theassump}{(A-$r$)}
\begin{assump}\label{assumAr}
The process $X$ satisfies Assumption~\ref{assumHr}, the
associated volatility process $c$ satisfies (H-$2$) and the processes
$b_t$ and, when $r\leq1$, $b'_t=b_t-\int\delta(t,z)
1_{\{\|\delta(t,z)\|\leq1\}} \la(dz)$ are c\`agl\`ad or c\`adl\`ag.
\end{assump}

The bigger $r$ is, the weaker Assumption~\ref{assumAr} is, and when
(A-$0$) holds the process $X$ has finitely many jumps on each finite
interval. Since we suppose in the theorems of the next section that $r<1$,
the last condition in (\ref{S-3})
implies that $b'_t$ is indeed well defined, and it is the ``genuine''
drift, in the sense that this is the drift after removing the sum
$\sum_{s\leq t}\De X_s$ of all jumps (which here are summable, and we even
have $\sum_{s\leq t}\|\De X_s\|^r<\infty$ a.s. here).

\section{The results}\label{sec-RES}

\subsection{A (seemingly) natural choice for the window $k_n$}

In order to define the estimators of the spot volatility, we need to
fix a
sequence $k_n$ of integers
and a sequence $u_n$ of cut-off levels in $(0,\infty]$. The
$\m_d^+$-valued variables $\Wc^n_i$ are defined, componentwise, as
%
\begin{equation}
\label{R-3} \widehat{c}{}^{n,lm}_i=\frac1{k_n
\De_n}\sum_{j=0}^{k_n-1}
\De^n_{i+j} X^l \De^n_{i+j}
X^m 1_{\{\|\De^n_{i+j}X\|\leq u_n\}},
\end{equation}
and they implicitly depend on $\De_n,k_n,u_n$.

A natural idea is to choose the sequence $k_n$ satisfying, as $n\to
\infty$,
%
\begin{equation}
\label{R-01} k_n\sim\frac{\te}{\rdn},\qquad \te\in(0,\infty).
\end{equation}
Indeed,\vspace*{1pt} one knows that $\widehat{c}{}^n_{[t/\De_n]}\toop c_t$ for any
$t$, as soon as
$k_n\to\infty$ and $k_n\De_n\to0$, and there is an
associated central limit theorem under Assumption~\ref{assumAr} for some $r<2$, with rate
$\min(1/\sqrt{k_n},1/\sqrt{k_n\De_n})$, which reaches its biggest value
$1/\De_n^{1/4}$ when\vadjust{\goodbreak}
$k_n\asymp1/\rdn$: this choice of $k_n$ ensures a balance between the involved
``statistical error'' which is of order $1/\sqrt{k_n}$, and the
variation of $c_t$ over the interval $[t,t+k_n\De_n]$, which
is of order $\sqrt{k_n\De_n}$ because $c_t$ is an It\^o semimartingale
(and even when it jumps); see~\cite{APPS,JP}.

By Theorem 9.4.1 of~\cite{JP}, and again as soon as
$k_n\to\infty$ and $k_n\De_n\to0$, one also knows that
%
\begin{equation}
\label{R-4} V(g)^n_t:=\De_n\sum
_{i=1}^{[t/\De_n]-k_n+1}g\bigl(\widehat{c}{}^n_i
\bigr)\quad\toucp\quad V(g)_t:=\int_0^tg(c_s)
\,ds
\end{equation}
(convergence in probability, uniform over each compact interval; by
convention $\sum_{i=a}^bv_i=0$ whenever $b<a$), as soon as the function
$g$ on $\m^+_d$ is continuous with $|g(x)|\leq K(1+\|x\|^p)$ for some
constants $K,p$, and under either one of the following three conditions:
%
\begin{eqnarray}
\label{R-5}
&&\begin{tabular}{p{310pt}}
$\bullet$
(A-0) holds, $X$ is continuous, $\frac{u_n}{\De
_n^{\ep}}\to
\infty$ for some $\ep<\frac12$ (e.g., $u_n\equiv\infty$);
\\
$\bullet$
(A-$r$) holds for some $r<2$ and $p\leq1$ and $u_n
\asymp\De_n^\vpi$ for some $\vpi\in(0,\frac12
)$;
\\
$\bullet$
(A-$r$) holds for some $r<2$ and $p>1$ and $u_n\asymp
\De_n^\vpi$ for some $\vpi\in[\frac{p-1}{2p-r},\frac12
)$.
\end{tabular}\hspace*{-26pt}
\end{eqnarray}
Notice the upper limit in definition (\ref{R-4}) of $V^n(g)_t$:
this is to ensure that $V^n(g)_t$ is actually computable from the
observations up to the time horizon~$t$. Note also that when $X$ is
continuous, the truncation in
(\ref{R-3}) is useless: one may use (\ref{R-3}) with $u_n\equiv
\infty$,
which reduces to (\ref{I-2}) in the one-dimensional case.

Now, we want to determine at which rate convergence (\ref{R-4})
takes place. This amounts to proving an associated central limit theorem.
Under the restriction $r<1$ and an appropriate choice of the truncation
levels, such a CLT is available for $V(g)^n$, with the rate $1/\rdn$,
but the
limit exhibits a bias term.

Below, $g$ is a smooth function on
$\m^+_d$, and the two first partial derivatives are denoted as
$\partial_{jk}g$ and $\partial^2_{jk,lm}g$, since any $x\in\m^+_d$ has
$d^2$ components $x^{jk}$. The family of all partial derivatives of
order $j$ is simply denoted as $\partial^jg$.

\begin{theo}\label{TR-01} Assume Assumption~\ref{assumAr} for some $r<1$. Let $g$ be a $C^3$
function on $\m^+_d$ such that
\[
\bigl\|\partial^jg(x)\bigr\|\leq K\bigl(1+\|x\|^{p-j}\bigr),\qquad
j=0,1,2,3,
\]
for some constants $K>0, p\geq3$. Either suppose that $X$ is continuous
and $u_n/\De_n^\ep\to\infty$ for some $\ep<1/2$ (e.g.,
$u_n\equiv\infty$, so there is no truncation at all), or suppose that
\[
u_n\asymp\De_n^\vpi,\qquad \frac{2p-1}{2(2p-r)}\leq
\vpi<\frac12.
\]
Then we have the finite-dimensional (in time) stable convergence in law
\[
\frac{1}{\rdn} \bigl(V(g)^n_t-V(g)_t
\bigr) \tolfs A^1_t+A^2_t+A^3_t+A^4_t+Z_t,
\]
where $Z$ is a process defined on an extension
$(\widetilde{\Omega},\widetilde{\f},(\widetilde{\f}_t)_{t\geq
0},\widetilde{\PP})$ of $(\Omega,\f,\break(\f_t)_{t\geq0},\PP)$, which
conditionally on $\f$ is a continuous
centered Gaussian martingale with variance
\[
\WE\bigl((Z_t)^2\mid\f\bigr) = \sum
_{j,k,l,m=1}^d\int_0^t
\partial_{jk} g(c_s) \,\partial_{lm}
g(c_s) \bigl(c_s^{jl}c_s^{km}+c_s^{jm}c_s^{kl}
\bigr)\,ds,
\]
and where
\begin{eqnarray*}
A^1_t&=&-\frac{\te}2 \bigl(g(c_0)+g(c_t)
\bigr),
\\
A^2_t&=&\frac1{2\te} \sum_{j,k,l,m=1}^d
\int_0^t\partial^2_{jk,lm}
g(c_s) \bigl(c^{jl}_sc^{km}_s+c^{jm}_sc^{kl}_s
\bigr)\,ds,
\\
A^3_t&=&-\frac{\te}{12} \sum
_{j,k,l,m=1}^d \int_0^t
\partial^2_{jk,lm} g(c_s) \Wc_s^{jk,lm}
\,ds,
\end{eqnarray*}
where $\Wc_s$ is the volatility process of
$c_t$,
\[
A^4_t=\te\sum_{s\leq t}G(c_{s-},
\De c_s)
\]
with $G(x,y)=\int_0^1
(g(x+wy)-(1-w)g(x)-wg(x+y))\,dw$.
\end{theo}

Note that $|G(x,y)|\leq K(1+\|x\|)^p \|y\|^2$, so the sum defining
$A_t^4$ is absolutely
convergent, and vanishes when $c_t$ is continuous.

The bias has four parts:
\begin{longlist}[(3)]
\item[(1)] The first part $A^1$ is a border effect, easily eliminated
by taking
%
\begin{equation}
\label{R-6} \WV(g)^n_t=V(g)^n_t+
\frac{(k_n-1)\De_n}2 \bigl(g\bigl(\widehat{c}{}^n_1\bigr)+g\bigl(
\widehat{c}{}^n_{[t/\De_n]-k_n+1} \bigr)\bigr)
\end{equation}
instead of $V(g)^n_t$: we then have
$\frac{1}{\rdn} (\WV(g)^n_t-V(g)_t)\tolfs A^2_t+A^3_t+A^4_t+Z_t$,
and this
convergence is even functional in time when $c_t$ is continuous.

\item[(2)] The second part $A^2$ is continuous in time and is present even
for the toy model $X_t=\sqrt{c} W_t$ with $c$ a constant and
$\De_n=\frac1n$ and $T=1$. In this simple case it can be interpreted as
follows: instead of taking the ``optimal'' $g(\widehat{c}_n)$ for
estimating $g(c)$,
with $\widehat{c}_n=\sum_{i=1}^n(\De^n_iX)^2$, one takes $\frac
1n\sum_{i=1}^n g(\widehat{c}{}^n_i)$
with $\widehat{c}{}^n_i$ a ``local'' estimator of $c$. This adds a
statistical error
which results in a bias.

\item[(3)] The third and fourth parts $A^3$ and $A^4$ are,
respectively,
continuous and purely discontinuous,\vadjust{\goodbreak} due to the continuous part and
to the jumps of the volatility process $c_t$ itself. These two biases
disappear if we take $\te=0$ in (\ref{R-01}) (with still $k_n\to
\infty$).
\end{longlist}

The only test function $g$ for which the biases $A^2,A^3,A^4$ disappear
is the identity $g(x)=x$. This is because, in this case, and up to
border terms, $\WV(g)^n_t$ is nothing but the realized quadratic
variation itself and the spot estimators $\widehat{c}{}^n_i$ actually
merge together
and disappear as such.

It is possible to consistently estimate $A^2_t,A^3_t,A^4_t$, and thus
de-bias $\WV(g)^n_t$ and obtain a CLT with a conditionally centered
Gaussian limit. Consistent estimators for $A^2_t$ are easy
to derive, since $A^2_t=V(f)_t$ for
the function $f(x)=\sum_{j,k,l,m}\partial^2_{jk,lm} g(x)
(x^{jl}x^{km}+x^{jm}x^{kl})$. Consistent estimators
for $A^3_t$ and $A^4_t$,
involving the volatility and the jumps of $c_t$, are more complicated to
describe, especially the last one, and also likely to have poor performances.
All the details about the way to remove the bias together with the proof
of Theorem~\ref{TR-01} can be found in~\cite{JR-tech}.

\subsection{A suitable window $k_n$}

In front of the difficulties involved in de-biasing the estimators $V(g)^n_t$
above, we in fact choose a window size $k_n$ smaller than the one in
(\ref{R-01}). Namely, we choose $k_n$ such that, as $n\to\infty$,
%
\begin{equation}
\label{R-1} k_n^3\De_n\to\infty,\qquad
k_n^2\De_n\to0.
\end{equation}
Of course, the second condition enables us to make the first and last
two bias
terms in Theorem~\ref{TR-01} vanish, which is technically very convenient.
However, it amplifies the first bias term, which becomes
the leading term in the difference $V(g)^n-V(g)$, and thus a prior de-biasing
is necessary if we want a rate $1/\rdn$. This leads us to consider the
following estimator:
%
\begin{eqnarray}
\label{R0-1} V'(g)^n_t&=&\De_n
\sum_{i=1}^{[t/\De_n]-k_n+1} \Biggl(g\bigl(
\widehat{c}{}^n_i\bigr)-\frac1{2k_n} \sum
_{j,k,l,m=1}^d \partial^2_{jk,lm}
g\bigl(\widehat{c}{}^n_i\bigr)\nonumber\\[-8pt]\\[-8pt]
&&\hspace*{161.7pt}{}\times \bigl(\widehat{c}{}^{n,jl}_i
\widehat{c}{}^{n,km}_i +\widehat{c}_i^{n,jm}
\widehat{c}_i^{n,kl} \bigr) \Biggr).\hspace*{-22pt}\nonumber
\end{eqnarray}

This estimator uses overlapping intervals, in the sense that we estimate
$c_{(i-1)\De_n}$ on the basis of the time window $((i-1)\De
_n,(i+k_n-1)\De_n]$,
and then sum over all $i$'s. Another version is indeed possible, which
does not use overlapping intervals and is as follows:
%
\begin{eqnarray}
\label{R0-101} V''(g)^n_t&=&k_n
\De_n\sum_{i=0}^{[t/k_n\De_n]-1} \Biggl(g
\bigl(\widehat{c}{}^n_{ik_n+1}\bigr)
\nonumber\\
&&\hspace*{73.3pt}{}-\frac1{2k_n} \sum_{j,k,l,m=1}^d
\partial^2_{jk,lm} g\bigl(\widehat{c}{}^n_{ik_n+1}
\bigr) \\
&&\hspace*{141pt}{}\times\bigl(\widehat{c}{}^{n,jl}_{ik_n+1} \widehat{c}{}^{n,km}_{ik_n+1}
+\widehat{c}_{ik_n+1}^{n,jm} \widehat{c}_{ik_n+1}^{n,kl}
\bigr) \Biggr).
\nonumber
\end{eqnarray}

We can now give the final version of our associated central
limit theorems.

\begin{theo}\label{TR-1} Assume Assumption~\ref{assumAr} for some $r<1$. Let $g$ be a $C^3$
function on $\m^+_d$ such that
%
\begin{equation}
\label{R-8} \bigl\|\partial^jg(x)\bigr\|\leq K\bigl(1+\|x\|^{p-j}
\bigr),\qquad j=0,1,2,3,
\end{equation}
for some constants $K>0, p\geq3$. Either suppose that $X$ is continuous
and $u_n/\De_n^\ep\to\infty$ for some $\ep<1/2$ (e.g.,
$u_n\equiv\infty$, so there is no truncation at all), or suppose that
%
\begin{equation}
\label{R-9} u_n\asymp\De_n^\vpi,\qquad
\frac{2p-1}{2(2p-r)}\leq\vpi<\frac12.
\end{equation}
Then under (\ref{R-1}) we have the two (functional in time) stable
convergences in law
%
\begin{equation}
\label{R-10}\quad \frac{1}{\rdn} \bigl(V'(g)^n-V(g)
\bigr) \tolls Z,\qquad \frac{1}{\rdn} \bigl(V''(g)^n-V(g)
\bigr) \tolls Z,
\end{equation}
where $Z$ is a process defined on an extension
$(\widetilde{\Omega},\widetilde{\f},(\widetilde{\f}_t)_{t\geq
0},\widetilde{\PP})$ of $(\Omega,\f,\break(\f_t)_{t\geq0},\PP)$, which
conditionally on $\f$ is a continuous
centered Gaussian martingale with variance
%
\begin{equation}
\label{R-11}\quad \WE\bigl((Z_t)^2\mid\f\bigr)=\sum
_{j,k,l,m=1}^d\int_0^t
\partial_{jk} g(c_s) \,\partial_{lm}
g(c_s) \bigl(c_s^{jl}c_s^{km}+c_s^{jm}c_s^{kl}
\bigr)\,ds.
\end{equation}
\end{theo}

\begin{rema}\label{RR-02} When $X$ jumps, the requirement (\ref{R-9})
is exactly the same as in Theorem~\ref{TR-01}, and it implies $r<1$. This
restriction is not a surprise, since
one needs $r\leq1$ in order to estimate the integrated
volatility by the (truncated) realized volatility, with a
rate of convergence $1/\rdn$. Indeed, it is shown in~\cite{JRe} that
if $r>1$, the optimal rate in the
minimax sense is $ (\rdn\operatorname{log}(1/\rdn) )^{-(2-r)/2}$.
When $r=1$ it
is likely that the CLT
still holds for an appropriate choice of the sequence~$u_n$, and with
another additional bias; see, for example,~\cite{V10} for a slightly different
context. Here we let this borderline case aside.
\end{rema}

\begin{rema}\label{RR-03} The limiting process $Z$ is the same
in both Theorems~\ref{TR-01} and~\ref{TR-1}, but in the latter case
the functional convergence always holds. It is also the same for
(the normalized versions of) the processes $V'(g)^n$ and $V''(g)^n$, which
is somewhat a surprise since in many instances using overlapping intervals
instead of nonoverlapping intervals results in a strictly smaller asymptotic
variance;
this is for example the case for multipower variations, see Theorem 11.2.1 in
\cite{JP}. However, in practice, it is probably advisable to use $V'(g)^n$
rather than $V''(g)^n$, because the former estimator is likely to be
less sensitive to way-off values of the spot estimators $\widehat
{c}{}^n_i$ than
the latter one, due to the ``smoothing'' embedded in its definition.
\end{rema}

\begin{rema}\label{RR-06} The $C^3$ property of $g$ is somewhat
restrictive, as, for example, in the one-dimensional case it rules out
the powers $g(x)=x^r$ with $r\in(0,3)\setminus\{1,2\}$. It could be proved
that, in the one-dimensional case again, and if the processes $c_t$
and $c_{t-}$ do not vanish (equivalently, the process $1/c_t$ is locally
bounded), the result still holds when $g$ is $C^3$ on $(0,\infty)$ and
satisfies (\ref{R-8}) with an arbitrary $p>0$: here again, the fact that
$1/c_t$ is locally bounded is also necessary for having a CLT for the
functionals of (\ref{I-1}) (say, with $k=1$) when the test function
$f$ is $C^1$ outside $0$ only.
\end{rema}

\begin{rema}\label{RR-07} One should compare this result with
those of Mykland and Zhang in~\cite{MZ}: in that paper [in which
only the continuous one-dimension\-al case and the test functions $g(x)=x^r$
are considered] the authors propose to take $k_n=k$ in (\ref{R-3}). Of course
(\ref{R-1}) fails, but $V(g)^n$ in this case is actually of the form
(\ref{I-1}) and a CLT holds for $\frac{1}{\rdn} (\al(g,k) V(g)^n-V(g))$
[without de-biasing term, but with an appropriate multiplicative factor
$\al(g,k)$, which is explicitly known]: the asymptotic variance
is bigger than in (\ref{R-11}), but approaches this value when $k$ is
large.

An advantage of Mykland--Zhang's approach is that when $g$ is positive, hence
$V(g)_t$ as well, the estimators are also positive. In contrast, $V'(g)^n_t$
in (\ref{R0-1}) may be negative even when $g\geq0$ everywhere. Thus if
this positivity issue is important for a specific application, taking
$k_n=k$ ``large'' and the estimator $\al(g,k) V(g)^n_t$ might be advisable,
although it seems to work only when $g$ is a power function.
Moreover, if $V'(g)^n_t$ is negative, it probably means that there is not
enough data in order to obtain a relevant estimation.
\end{rema}

It is simple to make this CLT ``feasible,'' that is, usable in practice
for determining a confidence interval for $V(g)_t$ at any time $t>0$.
Indeed, we can define the following function on $\m^+_d$:
%
\begin{equation}
\label{R-102} \Bh(x)=\sum_{j,k,l,m=1}^d
\partial_{jk} g(x) \,\partial_{lm} g(x) \bigl(x^{jl}x^{km}+x^{jm}x^{kl}
\bigr),
\end{equation}
which is continuous with $\Bh(x)\leq K(1+\|x\|^{2p-2})$, and nonnegative
(and positive at each $x$ such that $\partial g(x)\neq0$).
(\ref{R-9}) implies the last condition in (\ref{R-5}), and we have
$V(\Bh)^n\toucp V(\Bh)$, with $V(\Bh)_t$ being the right-hand
side of (\ref{R-11}). Then we readily deduce:

\begin{cor}\label{CR-1} Under the assumptions of the previous theorem,
for any $t>0$ we have the following stable convergence in law,
where $Y$ is an $\n(0,1)$ variable:
%
\begin{equation}
\label{R-25} \quad\frac{V'(g)^n_t-V(g)_t}{\sqrt{\De_n V(\Bh)^n_t}}\tols Y
\qquad\mbox{in restriction to the set $
\bigl\{V(\Bh)_t>0\bigr\}$},
\end{equation}
and the same holds with $V''(g)^n_t$ instead of $V'(g)^n_t$.
\end{cor}

\subsection{Optimality of the procedures}

We address now the question of the optimality of our procedures.

For simplicity, we restrict our attention to the one-dimensional case $d=1$.
We denote by $\s$ the class of all one-dimensional continuous semimartingales
$X$ of the form (\ref{I-3}), with $a,f$ being $C^3$ functions with bounded
derivatives with further $f$ bounded away from $0$, and $W,\BW$ being two
independent Brownian motions,
and $\Bb_t,\Bsi_t$ being Lebesgue square-integrable processes, optional
with respect to the filtration generated by $\BW$, and with $(\Bsi_t)^2$
bounded away from $0$. Such an $X$ satisfies
(A-$0$), with $\si_t=f(t,X_t,Y_t)$.

Let $t>0$. In the following, we say that a sequence of estimators
$(T^n_t)_{n\geq1}$ of $V(g)_t$ satisfy Property $\mathcal{P}$ over
$\s$ if:
\begin{longlist}[(2)]
\item[(1)] the estimator $T_t^n$
is a function of $(X_{i\De_n}\dvtx0\leq i\leq[t/\De_n])$;
\item[(2)] for any $X\in\s$, the variables
$\frac{1}{\rdn} (T^n_t-V(g)_t)$ converge stably in law to a limit
$Z'_t$ (depending
of $g$ of course), defined on an extension of the space.
\end{longlist}
The following theorem gives three small steps toward optimality.

\begin{theo}\label{TR-5} Let $d=1$ and $g$ be a $C^3$ function on $\R_+$
satisfying (\ref{R-8}) and which is strictly increasing, or
strictly decreasing.

\begin{longlist}[(a)]
\item[(a)]
For the parametric model $X_t=\si W_t$,
where $c_t=\si_t^2=c$ is a constant (the toy example of the \hyperref[sec-I]{Introduction}),
for any $t>0$, the estimators $V'(g)^n_t$ and $V''(g)^n_t$ are asymptotically
efficient (in Le Cam's sense) for estimating the number $tg(c)$.

\item[(b)] Let $(T^n_t)_{n\geq1}$ be a sequence of estimators satisfying
$\mathcal{P}$ over the class of continuous processes $X$ for which
\textup{(A-$0$)}
holds. Assume $Z'_t$ has a conditional variance of the form
%
\begin{equation}
\label{R-105} \WE\bigl(\bigl(Z'_t\bigr)^2
\mid\f\bigr)=\int_0^tH(c_s)\,ds
\end{equation}
for some nonnegative Borel function $H$. Then necessarily $H\geq\Bh$,
as given by (\ref{R-102}), and in particular,
%
\begin{equation}
\label{R-106} \WE\bigl(\bigl(Z'_t\bigr)^2
\mid\f\bigr)\geq\WE\bigl((Z_t)^2\mid\f\bigr).
\end{equation}

\item[(c)] The estimators $V'(g)^n_t$ and $V''(g)^n_t$
are optimal over $\s$ in the following sense: for any
sequence $(T^n_t)$ of estimators satisfying $\mathcal{P}$ over $\s$, the
limiting variable $Z'_t$ can be realized as $Z_t+Z_t''$, where $Z_t$ is the
limiting process in (\ref{R-10}), and the variable $Z''_t$ is
independent of
$Z_t$ conditionally on $\f$.
\end{longlist}
\end{theo}

Part (b) of Theorem~\ref{TR-5} shows in particular that the estimators
$U^n(f)_t$ given in (\ref{I-1}) for
estimating $g(x)=\E(f(\sqrt{x} U))$ have always an asymptotic
variance bigger than or equal to the variance (\ref{R-11}).

Part (c) states that our estimators achieve the lower bounds of Hajek
\textit{convolution theorem} over the class $\mathcal{S}$.
This convolution theorem for the subclass $\mathcal{S}$ is due to Cl\'ement,
Delattre and Gloter; see~\cite{CDG}. It in particular implies that for given
$t$, any rate optimal estimator over $\s$ has a limiting variance
which is
larger than those of $Z_t$ the limiting process in (\ref{R-10}).

So far, however, a ``general'' theory of optimality in our nonparametric
context seems still out of reach.

\begin{ex}[(Quarticity)]
Suppose $d=1$, and take $g(x)=x^2$, so
we want to estimate the quarticity $\int_0^tc_s^2\,ds$. In this case
an ``optimal'' estimator for the quarticity is
\[
\De_n \biggl(1-\frac2{k_n} \biggr)\sum
_{i=1}^{[t/\De_n]-k_n+1}\bigl(\widehat{c}{}^n_i
\bigr)^2.
\]
The asymptotic variance is $8\int_0^tc_s^4\,ds$, to be compared with the
asymptotic variance of the more usual estimators
$\frac1{3\De_n}\sum_{i=1}^{[t/\De_n]}(\De^n_iX)^4$, which is
$\frac{32}3\int_0^tc_s^4\,ds$.
\end{ex}

\begin{rema}\label{RR-15}
Although taking (\ref{R-1}) eliminates
the bias
terms $A^1_t$, $A^3_t$ and $A^4_t$ showing in Theorem~\ref{TR-01}, it might
be judicious to still eliminate the (asymptotically negligible) bias $A^1_t$
by adding to $V'(g)^n_t$ the same correction term
$\frac{(k_n-1)\De_n}2
(g(\widehat{c}{}^n_1)+g(\widehat{c}{}^n_{[t/\De_n]-k_n+1} )$ as
in (\ref{R-6}).\vspace*{2pt}

Due to their probable instability, it does not seem advisable, though, to
eliminate the biases $A^3_t$ and $A^4_t$ by using
(with the proper normalization) the method of~\cite{JR-tech}.
\end{rema}

\section{Proofs}\label{sec-P}

Under Assumption~\ref{assumAr}, not only do we have (\ref{S-1}), but we can write $c_t$ in
a similar fashion:
\begin{eqnarray*}
c_t&=&c_0+\int_0^t
\Wb_s\,ds+\int_0^t
\Wsi_s\,dW'_{s}+\int_0^t
\int_E \Wde(s,z) 1_{\{\|\Wde(s,z)\|\leq1\}} (\mu-\nu) (ds,dz)
\\
&&{} +\int_0^t\int
_E \Wde(s,z) 1_{\{\|\Wde(s,z)\|>1\}} \mu(ds,dz)
\end{eqnarray*}
(here, $W'$ is a $d^2$-dimensional Brownian motion, possibly correlated
with $W$).
Then, according to the localization Lemma 4.4.9 of~\cite{JP} [for the
assumption~(K) in that lemma], it is enough to show Theorem~\ref{TR-1}
under the following stronger assumption:

\renewcommand{\theassump}{(SA-$r$)}
\begin{assump}\label{assumSAr}
We have Assumption~\ref{assumAr}. Moreover we have, for a
$\la$-integrable function $J$ on $E$ and a constant $A$,
%
\begin{eqnarray}
\label{P-0}\qquad &\|b\|,\|\Wb\|,\|c\|,\|\Wc\|,J\leq A,\qquad \bigl\|\delta(\omega,t,z)
\bigr\|^r\leq J(z),&\nonumber\\[-8pt]\\[-8pt]
&\bigl\|\Wde(\omega,t,z)\bigr\|^2\leq J(z).&\nonumber
\end{eqnarray}
\end{assump}

In the sequel we suppose that $X$ satisfies Assumption~\ref{assumSAr}, and
also that (\ref{R-1}) holds: these assumptions are typically not recalled.
Below, all constants are denoted by $K$, and they vary from line to line.
They may implicitly depend on the process $X$ [usually through $A$ in
(\ref{P-0})]. When they depend on an additional parameter $p$, we
write $K_p$.

Recall the notation $b'_t$ in Assumption~\ref{assumAr}. We will usually replace
the discontinuous process $X$ by the continuous process
%
\begin{equation}
\label{P-3} X'_t=\int_0^tb'_s
\,ds+\int_0^t\si_s
\,dW_{s},
\end{equation}
connected with $X$ by $X_t=X_0+X'_t+\sum_{s\leq t}\De X_s$. Note that
$b'$ is
bounded, and without loss of generality we will use below its c\`adl\`ag
version.

\subsection{Estimates}\label{ssec-ARIS}

(1) First, we recall well-known estimates for $X'$
and $c$. Under (\ref{P-0}) and for $s,t\geq0$ and $q\geq0$, we have
%
\begin{eqnarray}
\label{P-251} \E\Bigl(\sup_{w\in[0,s]} \bigl\|X'_{t+w}-X'_t
\bigr\|^q\bigm|\f_t \Bigr)&\leq& K_q
s^{q/2},\nonumber\\
\bigl\| \E\bigl(X'_{t+s}-X'_t
\mid\f_s\bigr)\bigr\|&\leq& Ks,
\nonumber\\[-8pt]\\[-8pt]
\E\Bigl(\sup_{w\in[0,s]} \|c_{t+w}-c_t
\|^q\bigm|\f_t \Bigr) &\leq& K_q
s^{1\wedge(q/2)},\nonumber\\
\bigl\| \E(c_{t+s}-c_t\mid\f_s)
\bigr\|&\leq& Ks.
\nonumber
\end{eqnarray}

We need slightly more refined estimates for $X'$, and before giving them
we introduce some simplifying notation,
%
\begin{eqnarray}
\label{P-2} c^n_i&=&c_{(i-1)\De_n},\qquad
\f^n_i=\f_{(i-1)\De_n},
\nonumber\\
\eta_{t,s}&=&\sup\bigl(\bigl\|b'_{t+u}-b'_t
\bigr\|^2\dvtx  u\in[0,s] \bigr),\\
\eta^n_{i,j}&=&\sqrt
{\E\bigl(\eta_{(i-1)\De_n,j\De_n}\mid\f^{n}_i\bigr)},\qquad
\eta^n_i=\eta^n_{i,k_n}.
\nonumber
\end{eqnarray}

\begin{lem}\label{LP-6} We have
\begin{eqnarray*}\hspace*{-5pt}
&&\bigl|\E\bigl(\De^n_iX'^j
\De^n_iX'^m\mid
\f^{n}_i\bigr)-c^{n,jm}_i\De
_n \bigr| \\\hspace*{-5pt}
&&\qquad\leq K\De_n^{3/2}\bigl(\rdn+
\eta^n_{i,1}\bigr),
\\\hspace*{-5pt}
&&\bigl|\E\bigl(\De^n_iX'^j
\De^n_iX'^k
\De^n_iX'^l
\De^n_iX'^m \mid\f
^{n}_i \bigr)\,{-}\,\bigl(c^{n,jk}_ic^{n,lm}_i\,{+}\,c^{n,jl}_ic^{n,km}_i
\,{+}\,c^{n,jm}_ic^{n,kl}_i\bigr)
\De_n^2 \bigr|\\\hspace*{-5pt}
&&\qquad\leq K\De_n^{5/2}.
\end{eqnarray*}
\end{lem}

\begin{pf}
For simplicity we prove the result when $i=1$,
so $\De^n_1X'=X'_{\De_n}$, but upon shifting time the proof for $i>1$ is
the same.

First we have $X'_t=M_t+tb'_0+\int_0^t(b'_s-b'_0)\,ds$, where $M$ is a
martingale with $M_0=0$. Taking the $\f_0$-conditional expectation
thus yields
%
\begin{equation}
\label{P-503} \bigl\|\E\bigl(X'_t\mid\f_0
\bigr)-tb'_0 \bigr\|\leq t \eta_{0,t}.
\end{equation}
Next, It\^o's formula yields that $X'^j_tX'^m_t$ is the sum of a martingale
vanishing at $0$, plus
\begin{eqnarray*}
&&
b'^j_0 \int_0^t X_s'^m\,ds+b'^m_0 \int_0^t X_s'^j \,ds+\int_0^t
X_s'^m \bigl(b_s'^j-b_0'^j \bigr) \,ds\\
&&\qquad{} +\int_0^t X_s'^j
\bigl(b_s'^m-b_0'^m \bigr)\,ds +c^{jm}_0t+\int_0^t
\bigl(c^{jm}_s-c^{jm}_0\bigr)\,ds.
\end{eqnarray*}
Upon taking the conditional expectation, and using the Cauchy--Schwarz
inequality and the first and the last parts of (\ref{P-251}),
plus (\ref{P-503}), we readily deduce
%
\begin{equation}
\label{P-5030} \bigl|\E\bigl(X'^j_t
X'^m_t\mid\f_0
\bigr)-tc^{jm}_0 \bigr|\leq Kt^{3/2}(\sqrt{t}+
\eta_{0,t}).
\end{equation}

With $t=\De_n$, this gives the first claim. Finally, for
any indices $j_1,\ldots,j_4$ It\^o's formula yields a martingale $M$
vanishing at $0$ such that
%
\begin{eqnarray}
\label{P-504}\quad \prod_{l=1}^4
\De^n_1X'^{j_l}&=&M_{\De_n}
+\sum_{l=1}^p\int_0^{\De_n}b'^{j_l}_s
\prod_{1\leq m\leq p,m\neq l} X'^{j_m}_s
\,ds
\nonumber
\\
&&{}+\frac12 \sum_{1\leq l,l'\leq d,l\neq l'}c^{j_lj_{l'}}_0
\int_0^{\De_n}\prod_{1\leq m\leq4,m\neq l,l'}X'^{j_m}_s
\,ds
\\
&&{}+\frac12 \sum_{1\leq l,l'\leq d,l\neq l'}\int_0^{\De_n}
\bigl(c^{j_lj_{l'}}_s-c^{j_lj_{l'}}_0\bigr) \prod
_{1\leq m\leq4,m\neq l,l'} X'^{j_m}_s
\,ds.
\nonumber
\end{eqnarray}
Again, we take the $\f_0$-conditional expectation and we deal with the second,
the third and the last term in the right-hand side above by Fubini's
theorem and the Cauchy--Schwarz inequality. For the fourth term we use
(\ref{P-5030}), and a simple calculation yields the second claim.
\end{pf}

\begin{lem}\label{LP-1} For all $t>0$ we have
$\De_n\E(\sum_{i=1}^{[t/\De_n]}\eta^n_i )\to0$, and for
all $j,k$ such that
$j+k\leq k_n$ we have $\E(\eta^n_{i+j,k}\mid\f^{n}_i)\leq\eta^n_i$.
\end{lem}

\begin{pf}
The second claim follows from the definitions of
$\eta^n_i$ and $\eta^n_{i,j}$ and the Cauchy--Schwarz inequality.
For the first claim, we observe that $\mathbb{E}((\eta^n_{i})^2)$
is smaller than a constant always, and than
$\frac1{\Delta_{n}}\int_{(i-2)\Delta_{n}} ^{(i-1)\Delta_{n}
}\mathbb{E}((\eta_{s,2k_n+1})^2)\,ds$ when $i\geq2$. Hence by the
Cauchy--Schwarz inequality,
\begin{eqnarray*}
\Delta_{n}\mathbb{E} \Biggl(\sum_{i=1}^{[t/\Delta_{n}]}
\eta^n_{i} \Biggr) &\leq&\Biggl(t\mathbb{E} \Biggl(
\Delta_{n}\sum_{i=1}^{[t/\Delta_{n}]}\bigl(
\eta^n_{i}\bigr)^{2} \Biggr)
\Biggr)^{1/2} \\
&\leq&\biggl(Kt\De_n+\mathbb{E} \biggl(t\int
_{0}^{t}(\eta_{s,2k_n+1} )^2\,ds
\biggr) \biggr)^{1/2}.
\end{eqnarray*}
We have $\eta_{s,2k_n+1}\leq K$, and the c\`adl\`ag property of $b'$
yields that $\eta_{s,2k_n+1}(\omega)\to0$
for all $\omega$, and all $s$ except for
countably many strictly positive values (depending on $\omega$). Then, the
first claim follows by the dominated convergence theorem.
\end{pf}

(2) It is much easier (although unfeasible in practice) to
replace $\widehat{c}{}^n_i$ in (\ref{R-4}) by the estimators based on
the process
$X'$, as
given by (\ref{P-3}). Namely, we will replace $\widehat{c}{}^n_i$ by
the following:
\[
\widehat{c}'^n_i=\frac1{k_n
\De_n}\sum_{j=0}^{k_n-1}
\De^n_{i+j}X' \De^n_{i+j}X'^*.
\]

The comparison between $\widehat{c}{}^n_i$ and $\widehat{c}'^n_i$ is
based on the following
consequence of Lemma 13.2.6 of~\cite{JP}, applied with $F(x)=x x^*$, so
$k=1$ and $p'=s'=2$ and $s=1$ and $\ep=0$ (because $r<1$) with the notation
of that lemma. Namely, we have for all $q\geq1$ and for some sequence $a_n$
going to $0$,
\begin{eqnarray*}
&&
\E\bigl( \bigl\|\bigl(\De^n_iX \De^n_iX^*
\bigr) 1_{\{\|\De^n_iX\|\leq u_n\}}- \bigl(\De^n_iX'
\De^n_iX'^*\bigr) 1_{\{\|\De^n_iX'\|\leq u_n\}} \bigr\|
^q\bigr)  \\
&&\qquad\leq K_qa_n
\De_n^{(2q-r)\vpi+1}.
\end{eqnarray*}
Since $\E(\|\De^n_iX'\|^{2q})\leq K_q\De_n^q$ for any $q>0$ by classical
estimates, implying by Markov's inequality that
$\E(\|\De^n_iX'\|^{2q} 1_{\{\|\De^n_iX'\|>u_n\}})
\leq K\De_n^{q+q'(1-2\vpi)}$ for any $q'>0$, by taking $q'>\frac
1{1-2\vpi}$,
we then easily deduce
%
\begin{equation}
\label{P-6} \E\bigl(\bigl\|\widehat{c}{}^n_i-
\widehat{c}'^n_i\bigr\|^q \bigr)
\leq K_qa_n \De_n^{(2q-r)\vpi+1-q}.
\end{equation}

(3) Let us introduce the following $\R^d\otimes\R^d$-valued
variables:
%
\begin{eqnarray}
\label{P-20} \al^n_i&=&\De^n_iX'
\De^n_iX'^*-c^n_i
\De_n,\nonumber\\[-8pt]\\[-8pt]
\be^n_i&=&\widehat{c}'^n_i-c^n_i=
\frac1{k_n\De_n} \sum_{j=0}^{k_n-1}
\bigl(\al^n_{i+j}+\bigl(c^n_{i+j}-c^n_i
\bigr)\De_n \bigr).\nonumber
\end{eqnarray}

From (\ref{P-251}) we get that for all $q\geq0$,
%
\begin{equation}
\label{P-21} \E\bigl(\bigl\|\al^n_i\bigr\|^q\mid
\f^{n}_i\bigr)\leq K_q\De_n^q,\qquad
\bigl\|\E\bigl(\al^n_i\mid\f^{n}_i
\bigr)\bigr\|\leq K\De_n^{3/2}.
\end{equation}
This and the Burkholder--Gundy and H\"older inequalities
give us, for $q\geq2$, that
$\E(\|{\sum_{j=0}^{k_n-1}\al^n_{i+j}}\|^q\mid\f^n_i )
\leq K_q\De_n^qk_n^{q/2}$. This\vspace*{2pt} and (\ref{P-251}) and again H\"older's
inequality yield
%
\begin{equation}
\label{P-23} q\geq2\quad\Rightarrow\quad\E\bigl(\bigl\|\be^n_i
\bigr\|^q\mid\f^{n}_i\bigr)\leq K_q
\bigl(k_n^{-q/2}+k_n\De_n \bigr).
\end{equation}
Lemma~\ref{LP-6} allows us for better estimates for $\al^n_i$, namely
%
\begin{eqnarray}
\label{P-22} &\bigl\| \E\bigl(\al^n_i\mid\f^{n}_i
\bigr) \bigr\| \leq K\De_n^{3/2} \bigl(\rdn+\eta^n_{i,1}
\bigr),&
\nonumber\\[-8pt]\\[-8pt]
&\bigl|\E\bigl(\al^{n,jk}_i\al_i^{n,lm}\mid
\f^{n}_i\bigr) -\bigl(c^{n,jl}_ic^{n,km}_i+c^{n,jm}_ic^{n,kl}_i
\bigr)\De_n^2 \bigr|\leq K\De_n^{5/2}.&
\nonumber
\end{eqnarray}

\begin{lem}\label{LP-7} We have
\begin{eqnarray*}
&&\bigl\| \E\bigl(\be^n_i\mid\f^{n}_i
\bigr) \bigr\| \leq K \rdn\bigl(k_n\rdn+\eta^n_i
\bigr),
\\
&&\biggl|\E\bigl(\be^{n,jk}_i \be_i^{n,lm}\mid
\f^{n}_i\bigr) -\frac1{k_n}
\bigl(c^{n,jl}_ic^{n,km}_i+c^{n,jm}_ic^{n,kl}_i
\bigr) \biggr| \\
&&\qquad\leq K\rdn\bigl(k_n^{-1/2}+k_n\rdn+
\eta^n_i\bigr).
\end{eqnarray*}
\end{lem}

\begin{pf}
The first claim follows from (\ref{P-251}), (\ref{P-22})
and the last part of Lem\-ma~\ref{LP-1}. For the second one, we set
$\xi^n_i=c^{n,jl}_ic^{n,km}_i+c^{n,jm}_ic^{n,kl}_i$ and
$\ze^n_{i,j}=\al^n_{i+j}+(c^n_{i+j}-c_i^n)\De_n$ and write
$\be^{n,jk}_i\be_i^{n,lm}$ as
%
\begin{eqnarray}
\label{P-2201}
&&\frac1{k_n^2\De_n^2}
\sum_{u=0}^{k_n-1}\ze_{i,u}^{n,jk}
\ze^{n,lm}_{i,u} +\frac1{k_n^2
\De_n^2}\sum_{u=0}^{k_n-2}
\sum_{v=u+1}^{k_n-1} \ze^{n,jk}_{i,u}
\ze^{n,lm}_{i,v} \nonumber\\[-8pt]\\[-8pt]
&&\qquad{}+\frac1{k_n^2
\De_n^2}\sum_{u=0}^{k_n-2}
\sum_{v=u+1}^{k_n-1} \ze^{n,lm}_{i,u}
\ze^{n,jk}_{i,v}.\nonumber
\end{eqnarray}

First, we have
\[
\bigl|\ze_{i,u}^{n,jk}\ze^{n,lm}_{i,u}-
\al_{i+u}^{n,jk}\al^{n,lm}_{i+u}\bigr| \leq2
\De_n\bigl\|c^n_{i+u}-c^n_i
\bigr\| \bigl\|\al^n_{i+u}\bigr\|+ \De_n^2
\bigl\|c^n_{i+u}-c^n_i\bigr\|^2,
\]
whose $\f^n_i$-conditional expectation is less than
$K\De_n^{5/2}k_n^{1/2}$ by (\ref{P-251}) and (\ref{P-21}).
The boundedness of $c_t$ and (\ref{P-251})
yield $ |\E(\xi_{i+u}^n\mid\f^{n}_i)-\xi^n_i |
\leq K k_n\De_n$. Then (\ref{P-22}) gives us that the $\f^n_i$-conditional
expectation of the first term in (\ref{P-2201}), minus $\frac1{k_n}
\xi^n_i$,
is less than $K\rdn/\sqrt{k_n}$.

Second, (\ref{P-251}) and (\ref{P-22}), plus the first claim of
Lemma~\ref{LP-6}, yield, when $0\leq u<v<k_n$,
\begin{eqnarray*}
\bigl|\E\bigl(\ze^{n,jk}_{i,v}\mid\f^n_{i+u+1}
\bigr)-\bigl(c^{n,jk}_{i+u+1}-c^{n,jk}_i
\bigr)\De_n\bigr| &\leq& K\De_n^{3/2}
\bigl(k_n\rdn+\eta^n_{i+v,1}\bigr),
\\
\bigl|\E\bigl(\al^{n,lm}_{i+u} \bigl(c^{n,jk}_{i+u+1}-c^{n,jk}_{i+u}
\bigr)\mid\f^n_{i+u}\bigr)\bigr| &\leq& K\De_n^{3/2}
\bigl(\rdn+\eta^n_{i+u,1}\bigr),
\\
\bigl|\E\bigl(\al^{n,lm}_{i+u} \bigl(c^{n,jk}_{i+u}-c^{n,jk}_i
\bigr)\mid\f^n_{i+u}\bigr)\bigr| &\leq& K\De_n^{3/2}
\bigl(\rdn+\eta^n_{i+u,1}\bigr),
\\
\bigl|\E\bigl(\bigl(c^{n,lm}_{i+u}-c^{n,lm}_i
\bigr) \bigl(c^{n,jk}_{i+u+1}-c^{n,jk}_i
\bigr)\mid\f^n_{i}\bigr)\bigr| &\leq& Kk_n
\De_n.
\end{eqnarray*}
Since $\ze^n_{i+u}$ is $\f^n_{i+u+1}$-measurable, and using (\ref
{P-21}) and
the second part of Lem\-ma~\ref{LP-1}, the $\f^{n}_i$-conditional
expectation of
the last term of (\ref{P-2201}) is smaller than $K \rdn(k_n\rdn
+\eta^n_i)$.
The same is obviously true for the second term, and
we readily deduce the second claim of the lemma.
\end{pf}

\subsection{\texorpdfstring{Proof of Theorem \protect\ref{TR-1}}{Proof of Theorem 3.2}}\label{ssec-2}

Using the key property $\widehat{c}'^n_i=c^n_i+\be^n_i$ and the definition
(\ref{P-20}) of $\be^n_i$, a simple calculation shows the decomposition
$\frac{1}{\rdn} (V'(g)^n_t-V(g)_t)=\sum_{j=1}^5V^{n,j}_t$,
as soon as $t>k_n\De_n$, and where
\begin{eqnarray*}
V^{n,1}_t&=&\rdn\sum_{i=1}^{[t/\De_n]-k_n+1}
\Biggl(g\bigl(\widehat{c}{}^n_i\bigr)-g\bigl(
\widehat{c}'^n_i\bigr)
\\
&&\hspace*{77.5pt}{}-\frac1{2k_n}\sum_{j,k,l,m=1}^d
\bigl(\partial^2_{jk,lm} g\bigl(\widehat{c}{}^n_i
\bigr) \bigl(\widehat{c}{}^{n,jl}_i \widehat{c}{}^{n,km}_i
+\widehat{c}_i^{n,jm} \widehat{c}_i^{n,kl}
\bigr)
\\
&&\hspace*{128.5pt}{}-\partial^2_{jk,lm} g\bigl(\widehat{c}'^n_i
\bigr) \bigl(\widehat{c}'^{n,jl}_i
\widehat{c}'^{n,km}_i+\widehat{c}_i'^{n,jm}
\widehat{c}_i'^{n,kl}\bigr) \bigr) \Biggr),
\\
V^{n,2}_t&=&\frac{1}{\rdn} \sum
_{i=1}^{[t/\De_n]-k_n+1} \int_{(i-1)\De_n}^{i\De_n}
\bigl(g\bigl(c_i^n\bigr)-g(c_s)\bigr)\,ds\\
&&{} -
\frac{1}{\rdn} \int_{\De_n([t/\De_n]-k_n+1)}^tg(c_s)
\,ds,
\\
V^{n,3}_t&=&\rdn\sum_{i=1}^{[t/\De_n]-k_n+1}
\sum_{l,m=1}^d \partial_{lm}g
\bigl(c^n_i\bigr) \frac{1}{k_n} \sum
_{u=0}^{k_n-1}\bigl(c_{i+u}^{n,lm}-c_i^{n,lm}
\bigr),
\\
V^{n,4}_t&=&\rdn\sum_{i=1}^{[t/\De_n]-k_n+1}
\Biggl(g\bigl(c^n_i+\be^n_i
\bigr)-g\bigl(c^n_i\bigr) -\sum
_{l,m=1}^d \partial_{lm}g
\bigl(c^n_i\bigr) \be^{n,lm}_i
\\
&&\hspace*{77.5pt}{}-\frac1{2k_n}\sum_{j,k,l,m=1}^d
\partial^2_{jk,lm} g\bigl(c^n_i+
\be^n_i\bigr)\\
&&\hspace*{145.6pt}{}\times \bigl(\bigl(c^{n,jl}_i+
\be^{n,jl}_i\bigr) \bigl(c^{n,km}_i+
\be^{n,km}_i\bigr)
\\
&&\hspace*{162.3pt}{}+\bigl(c^{n,jm}_i+\be^{n,jm}_i
\bigr) \bigl(c^{n,kl}_i+\be^{n,kl}_i
\bigr) \bigr) \Biggr),
\\
V^{n,5}_t&=&\frac1{k_n\rdn} \sum
_{i=1}^{[t/\De_n]-k_n+1} \sum_{l,m=1}^d
\partial_{lm}g\bigl(c^n_i\bigr) \sum
_{u=0}^{k_n-1}\al^{n,lm}_{i+u}.
\end{eqnarray*}

The leading term is $V^{n,5}$, and the first claim in (\ref{R-10}), about
$V'(g)^n$, is a consequence of the following two lemmas:

\begin{lem}\label{LP-10} For $v=1,2,3,4$ we have $V^{n,v}
\stackrel{\mathit{u.c.p.}}{\Longrightarrow}0$.
\end{lem}

\begin{lem}\label{LP-11} With $Z$ as in Theorem~\ref{TR-1}, we have the
functional stable convergence in law
%
\begin{equation}
\label{P-50} V^{n,5}\tolls Z.
\end{equation}
\end{lem}

\begin{pf*}{Proof of Lemma~\ref{LP-10}}
\textit{The case} $v=1$: We define functions
$h_n$ on $\m^+_d$ by
\[
h_n(x)=g(x)-\frac1{2k_n}\sum
_{j,k,l,m=1}^d \partial^2_{jk,lm} g(x)
\bigl(x^{jl}x^{km}+x^{jm}x^{kl}\bigr).
\]
From (\ref{R-8}) we obtain $|h_n(x)-h_n(y)|\leq K(1+\|y\|)^{p-1})\|
x-y\|
+K\|x-y\|^p$ (uniformly in $n$). So if $\eta^n_i$ is the $i$th summand
in the
definition of $V^{n,1}_t$, we get
\[
\bigl|\eta^n_i\bigr|\leq K\bigl(1+\bigl\|\widehat{c}{}^n_i
\bigr\|^{p-1}+\bigl\|\widehat{c}'^n_i\bigr\|
^{p-1}\bigr)\bigl\|\widehat{c}{}^n_i-
\widehat{c}'^n_i\bigr\| +K\bigl\|
\widehat{c}{}^n_-\widehat{c}'^n_i
\bigr\|^p.
\]
Recalling the last part of (\ref{P-21}), and by (\ref{P-6}), H\"older's
inequality and the fact that $(2p-r)\vpi+1-p<\frac1q ((2q-r)\vpi+1-q)$
when $q>1$ is small enough, because $\vpi<\frac12$, we deduce
$\E(|g(\widehat {c}{}^n_i)-g(\widehat{c}'^n_i)|)\leq
Ka_n\De_n^{(2p-r)\vpi+1-p}$ and thus
\[
\E\Bigl(\sup_{s\leq t} \bigl|V^{n,1}_s\bigr| \Bigr)
\leq Kta_n\De_n^{(2p-r)\vpi+1/2-p}.
\]
In view of (\ref{R-9}), we deduce the result for $v=1$.

\textit{The case} $v=2$: Since $g(c_s)$ is bounded, it is obvious that
the absolute value of the last term in $V^{n,2}_t$ is smaller than
$Kk_n\rdn$,
which goes to $0$ by (\ref{R-1}). Since $g$ is~$C^2$, the convergence
of the first term in $V^{n,2}_t$ to $0$ in
probability, locally uniformly in $t$, is well known; see, for example,
the proof of
(5.3.24) in~\cite{JP}, in which one replaces $\rho_{c_s}(f)$ by $g(c_s)$.
Thus the result holds for $v=2$.

\textit{The case} $v=3$: Letting $\ze^n_i=\sum_{l,m=1}^d
\partial_{lm}g(c^n_i) \frac{1}{k_n}
\sum_{u=0}^{k_n-1}(c_{i+u}^{n,lm}-c_i^{n,lm})$ be the $i$th summand in the
definition\vadjust{\goodbreak} of $V^{n,3}_t$, and $N(n,j,t)$ be the
integer part of $([t/\De_n]-k_n-j+1)/k_n$, we have
\[
V^{n,3}_t=\rdn\sum_{j=1}^{k_n}H(j)^n_t\qquad
\mbox{where } H(j)^n_t=\sum_{i=0}^{N(n,j,t)}
\ze^n_{j+k_ni}.
\]
From (\ref{P-251}) and the Cauchy--Schwarz inequality, we get
\[
\bigl|\E\bigl(\ze^n_i\mid\f^n_i
\bigr)\bigr|\leq Kk_n\De_n, \qquad\E\bigl(\bigl|\ze^n_i\bigr|^2
\mid\f^n_i\bigr)\leq Kk_n\De_n.
\]
Then Doob's inequality, the $\f^n_{j+k_n(i+1)}$-measurability of
$\ze^n_{j+k_ni}$, and $N(n,j,t)\leq t/k_n\De_n$ imply
\begin{eqnarray*}
\E\Bigl(\sup_{s\leq t}\bigl|H(j)^n_s\bigr| \Bigr)
&\leq&\sum_{i=0}^{N(n,j,t)}\E\bigl(\bigl|\E\bigl(
\ze^n_{j+k_ni}\mid\f^n_{j+k_ni}\bigr)\bigr|
\bigr)\\
&&{}+ \Biggl(4\sum_{i=0}^{N(n,j,t)}\E\bigl(
\bigl(\ze(j)^n_{j+k_ni}\bigr)^2\bigr)
\Biggr)^{1/2}\\
&\leq& K(t+\sqrt{t} ).
\end{eqnarray*}
Since $|V^{n,3}_t|\leq\rdn\sum_{j=1}^{k_n}|H(j)^n_t|$ and $k_n\rdn
\to0$, we
deduce the result for $v=3$.

\textit{The case} $v=4$: The $i$th summand in the definition of $V^{n,4}_t$
is $v^n_i+w^n_i$, where
\begin{eqnarray*}
v^n_i&=&\frac12\sum_{j,k,l,m=1}^d
\partial^2_{jk,lm} g\bigl(c^n_i\bigr)
\biggl(\be^{n,jk}_i \be^{n,lm}_i-
\frac1{k_n} \bigl(c^{n,jl}_i c^{n,km}_i
+c_i^{n,jm} c_i^{n,kl} \bigr) \biggr),
\\
\bigl|w^n_i\bigr|&\leq& K\bigl(1+\bigl\|\be^n_i
\bigr\|^{p-3}\bigr)\bigl\|\be^n_i\bigr\|^3+
\frac{K}{k_n} \bigl(1+\bigl\|\be^n_i\bigr\|^{p-1}
\bigr) \bigl\|\be^n_i\bigr\|
\end{eqnarray*}
[use (\ref{R-8}) and $\|c_t\|\leq K$ repeatedly], and we thus have
$V^{n,4}_t=G^n_t+\sum_{j=1}^{k_n}H(j)^n_t$, with $N(n,j,t)$ as in the
previous step and
\begin{eqnarray*}
G^n_t&=&\rdn\sum_{i=1}^{[t/\De_n]-k_n+1}
\bigl(w^n_i+\E\bigl(v^n_i\mid\f
^n_i\bigr)\bigr),
\\
H(j)^n_t&=&\sum_{i=0}^{N(n,j,t)}
\ze(j)^n_i,\qquad \ze(j)^n_i=\rdn
\bigl(v^n_{j+k_ni}-\E\bigl(v^n_{j+k_ni}\mid
\f^n_{j+k_ni}\bigr) \bigr).
\end{eqnarray*}

In view of Lemma~\ref{LP-7} and (\ref{P-23}), plus
H\"older's inequality, we have
\begin{eqnarray*}
\bigl|\E\bigl(v^n_i\mid\f^n_i
\bigr)\bigr|&\leq& K\rdn\bigl(k_n\rdn+\eta^n_i
\bigr),\\
\E\bigl(\bigl|w^n_i\bigr|\bigr)&\leq& K \biggl(
\frac1{k_n^{3/2}}+k_n\De_n+
\frac{\rdn
}{\sqrt{k_n}} \biggr),
\end{eqnarray*}
and thus (\ref{R-1}) and Lemma~\ref{LP-1} yield
\[
\E\Bigl(\sup_{s\leq t}\bigl|G^n_s\bigr| \Bigr)
\leq\E\Biggl(\sum_{i=1}^{[t/\De_n]}\rdn
\bigl(\bigl|w^n_i\bigr|+\bigl|\E\bigl(v^n_i\mid
\f^n_i\bigr)\bigr| \bigr) \Biggr)\to0.
\]
Moreover (\ref{P-23}) and $k_n^{-2}\leq Kk_n\De_n$ yield
$\E(|\ze(j)^n_i|^2)\leq K\De_n^2k_n$, whereas
$\ze(j)^n_i$ is a martingale increment for the
filtration $(\f^n_{j+k_ni})_{i\geq0}$, hence Doob's inequality and
$N(n,j,t)\leq t/k_n\De_n$ imply
\[
\E\Bigl(\sup_{s\leq t}\bigl|H(j)^n_s\bigr| \Bigr)
\leq\Biggl(\sum_{i=0}^{N(n,j,t)}\E\bigl(\bigl(
\ze(j)^n_i\bigr)^2\bigr)
\Biggr)^{1/2}\leq Kt \De_n.
\]
Since $|V^{n,4}_t|\leq|G^n_t|+\sum_{j=1}^{k_n}|H(j)^n_t|$, we deduce
the result for $v=4$.
\end{pf*}

\begin{pf*}{Proof of Lemma~\ref{LP-11}}
We can rewrite $V^{n,5}$ as
\[
V^{n,5}_t=\frac{1}{\rdn}\sum
_{i=1}^{[t/\De_n]}\sum_{l,m=1}^d
w^{n,lm}_i \al^{n,lm}_i,
\]
where
\[
w^{n,lm}_i=\frac1{k_n}\sum
_{j=(i-[t/\De_n]+k_n-1)^+}^{(i-1)\wedge(k_n-1)} \partial_{lm} g
\bigl(c_{i-j}^n\bigr).
\]
Observe that $w^n_i$ and
$\al^n_i$ are measurable with respect to $\f^n_i$ and $\f^n_{i+1}$,
respectively, so by Theorem IX.7.28 of~\cite{JS} (with $G=0$ and $Z=0$ in
the notation of that theorem) it suffices to prove the
following four convergences in
probability, for all $t>0$ and all component indices:
%
\begin{eqnarray}
\label{P-51}
&&\hspace*{23.8pt}\frac{1}{\rdn}\sum_{i=1}^{[t/\De_n]-k_n+1}w^{n,lm}_i
\E\bigl(\al^{n,lm}_i\mid\f^n_i
\bigr)\toop 0,
\\
%
\label{P-52}
&&\frac1{\De_n} \sum_{i=1}^{[t/\De_n]-k_n+1}
w^{n,jk}_i w^{n,lm}_i \E\bigl(
\al^{n,jk}_i \al^{n,lm}_i\mid
\f^n_i\bigr) \nonumber\\[-8pt]\\[-8pt]
&&\qquad\toop\int_0^t
\partial_{jk} g(c_s) \,\partial_{lm}
g(c_s) \bigl(c_s^{jl}c_s^{km}+c_s^{jm}c_s^{kl}
\bigr)\,ds,\nonumber
\\
%
\label{P-53}
&&\hspace*{28.1pt}\frac1{\De_n^2}\sum
_{i=1}^{[t/\De_n]-k_n+1}\bigl\|w^n_i
\bigr\|^4 \E\bigl(\bigl\|\al^n_i\bigr\|^4\mid
\f^n_i\bigr) \toop 0,
\\
%
\label{P-54}
&&\frac{1}{\rdn}\sum_{i=1}^{[t/\De_n]-k_n+1}w^{n,lm}_i
\E\bigl(\al^{n,lm}_i \De^n_iN\mid
\f^n_i\bigr) \toop 0,
\end{eqnarray}
where $N=W^k$ for some $k$, or is an arbitrary bounded
martingale, orthogonal to~$W$.

Lemma~\ref{LP-1}, (\ref{P-21}), (\ref{P-22}) and the property
$\|w^n_i\|\leq K$ readily
imply (\ref{P-51}) and (\ref{P-53}). In view of the form of $\al^n_i$,
a usual argument (see, e.g.,~\cite{JP}) shows that in fact
$\E(\al^{n,lm}_i \De^n_iN\mid\f^n_i)=0$ for all $N$ as above, and hence
(\ref{P-54}) holds.\looseness=-1

For (\ref{P-52}), by (\ref{P-22}) it suffices to prove that
\begin{eqnarray*}
&&
\De_n \sum_{i=1}^{[t/\De_n]-k_n+1}
w^{n,jk}_i w^{n,lm}_i
\bigl(c^{n,jl}_ic^{n,km}_i+c^{n,jm}_ic^{n,kl}_i
\bigr) \\
&&\qquad\toop\int_0^t \partial_{jk}
g(c_s) \,\partial_{lm} g(c_s)
\bigl(c_s^{jl}c_s^{km}+c_s^{jm}c_s^{kl}
\bigr)\,ds.
\end{eqnarray*}
In view of the definition of $w^n_i$, for each $t$ we have
$w^{n,jk}_{i(n,t)}\to\partial_{jk} g(c_t)$ and $c^{n,jk}_{i(n,t)}\to
c^{jk}_t$ almost surely if $|i(n,t)\De_n-t|\leq k_n\De_n$, and the above
convergence follows by the
dominated convergence theorem, thus ending the proof of (\ref
{P-50}).
\end{pf*}

\begin{pf*}{Proof of the second claim in (\ref{R-10})}
The proof is basically
the same as for the first claim. We have the decomposition
$\frac{1}{\rdn} (V''(g)^n_t-V(g)_t)=\sum_{j=1}^5\BV^{n,j}_t$, where
\begin{eqnarray*}
\BV^{n,1}_t&=&k_n\rdn\sum
_{i=0}^{[t/k_n\De_n]-1} \bigl(g\bigl(\widehat{c}{}^n_{k_ni+1}
\bigr)-g\bigl(\widehat{c}'^n_{k_ni+1}\bigr)
\bigr),
\\
\BV^{n,2}_t&=&\frac{1}{\rdn} \sum
_{i=0}^{[t/k_n\De_n]-1} \int_{k_ni\De_n}^{k_n(i+1)\De_n}
\bigl(g(c_{k_ni\De_n})-g(c_s)\bigr)\,ds \\
&&{}-\frac{1}{\rdn} \int
_{k_n\De_n([t/k_n\De_n])}^tg(c_s)\,ds,
\\
\BV^{n,3}_t&=&k_n\rdn\sum
_{i=0}^{[t/k_n\De_n]-1} \sum_{l,m=1}^d
\partial_{lm}g\bigl(c^n_{k_ni+1}\bigr)
\frac{1}{k_n} \sum_{u=0}^{k_n-1}
\bigl(c_{k_ni+1+u}^{n,lm}-c_{k_ni+1}^{n,lm}\bigr),
\\
\BV^{n,4}_t&=&k_n\rdn\\
&&{}\times\sum
_{i=0}^{[t/k_n\De_n]-1} \Biggl(g\bigl(c^n_{k_ni+1}+
\be^n_{k_ni+1}\bigr)-g\bigl(c^n_{k_ni+1}
\bigr) \\
&&\hspace*{82.2pt}\hspace*{-22.5pt}{}-\sum_{l,m=1}^d \partial_{lm}g
\bigl(c^n_{k_ni+1}\bigr) \be^{n,lm}_{k_ni+1}
\\
&&\hspace*{82.2pt}\hspace*{-22.5pt}{}-\frac1{2k_n}\sum_{j,k,l,m=1}^d
\partial^2_{jk,lm} g\bigl(c^n_{k_ni+1}+
\be^n_{k_ni+1}\bigr) \\
&&\hspace*{149.8pt}\hspace*{-22.5pt}{}\times\bigl(\bigl(c^{n,jl}_{k_ni+1}+
\be^{n,jl}_{k_ni+1}\bigr) \bigl(c^{n,km}_{k_ni+1}+
\be^{n,km}_{k_ni+1}\bigr)
\\
&&\hspace*{143.1pt}{} +\bigl(c^{n,jm}_{k_ni+1}+\be^{n,jm}_{k_ni+1}
\bigr) \bigl(c^{n,kl}_{k_ni+1}+\be^{n,kl}_{k_ni+1}
\bigr) \bigr) \Biggr),
\\
\BV^{n,5}_t&=&\frac1{\rdn} \sum_{i=0}^{[t/k_n\De_n]-1}
\sum_{l,m=1}^d \partial_{lm}g
\bigl(c^n_{k_ni+1}\bigr) \sum_{u=0}^{k_n-1}
\al^{n,lm}_{k_ni+u+1}.
\end{eqnarray*}
The proofs of Lemmas~\ref{LP-10} and~\ref{LP-11} carry over to $\BV^{n,v}$
instead of $V^{n,v}$, for $v=1,2,3,4,5$, almost word for word, except for
the following points:

(1) For Lemma~\ref{LP-10}, cases $v=3,4$, there is no need to consider
the $k_n$ processes $H(j)^n$; a single process $H^n$ is enough, and the proof
is simpler.

(2) For Lemma~\ref{LP-10}, case $v=2$, the proof of the u.c.p. convergence
to $0$ of the first term in the definition of $\BV^{n,2}$ should be
reworked as follows: the $i$th summand $\ze^n_i$ in this term is
$\f^n_{k_n(i+1)}$-measurable, and by (\ref{R-8}) and (\ref{P-251}) it
satisfies
\[
\bigl|\E\bigl(\ze^n_i\mid\f^n_{k_ni}
\bigr)\bigr|\leq K(k_n\De_n)^2,\qquad \E\bigl(\bigl|
\ze^n_i\bigr|^2\mid\f^n_{k_ni}
\bigr)\leq K(k_n\De_n)^3.
\]
Then the claim follows from the usual
martingale argument and $k_n\rdn\to0$.

(3) For Lemma~\ref{LP-11}, we have
\[
\BV^{n,5}_t=\frac{1}{\rdn}\sum
_{i=1}^{k_n[t/k_n\De_n]} \sum_{l,m=1}^d
\partial_{lm} g\bigl(c_{1+k_n[(j-1)/k_n]}^n\bigr)
\al^{n,lm}_i,
\]
and the rest of the proof is similar.
\end{pf*}

\subsection{\texorpdfstring{Proof of Theorem \protect\ref{TR-5}\textup{(a)} and \textup{(b)}}
{Proof of Theorem 3.8(a) and (b)}}
(a) is almost obvious: indeed, $\BV(g)^n_t$ converges with the rate
$\frac{1}{\rdn}$ and is asymptotically normal with asymptotic variance
$2tg'(c)^2c^2$ ($g'$ is the derivative of $g$). However, since $g$ is
one-to-one, the model index by the new parameter $tg(c)$ is regular,
and the MLE is $tg(\widehat{c}_n)$, where $\widehat{c}_n=\sum
_{i=1}^{[t/\De_n]}(\De^n_i
X)^2$, and clearly
$tg(\widehat{c}_n)$ has the same asymptotic properties as $\BV
(g)^n_t$: this
proves the result.

(b) is also obvious: the properties of $T^n_t$ hold for all continuous
processes $X$ satisfying (A-$0$). Then,\vspace*{1pt} using the toy model of (a), the
optimality proved above implies that $tH(c)\geq2tg'(c)^2c^2$ for any
constant $c>0$, that is, $H\geq\Bh$.

Finally, (c) is exactly Theorem 3 of~\cite{CDG} applied to the present
setting.

\section*{Acknowledgments}

We thank two anonymous referees for very constructive comments, which
led to substantial improvements. We also thank Jia Li for pointing out
some problems in an earlier draft.



\printaddresses


\begin{thebibliography}{13}

\bibitem{APPS}
\begin{barticle}[mr]
\bauthor{\bsnm{Alvarez},~\bfnm{Alexander}\binits{A.}},
  \bauthor{\bsnm{Panloup},~\bfnm{Fabien}\binits{F.}},
  \bauthor{\bsnm{Pontier},~\bfnm{Monique}\binits{M.}} \AND
  \bauthor{\bsnm{Savy},~\bfnm{Nicolas}\binits{N.}}
(\byear{2012}).
\btitle{Estimation of the instantaneous volatility}.
\bjournal{Stat. Inference Stoch. Process.}
\bvolume{15}
\bpages{27--59}.
\bid{doi={10.1007/s11203-011-9062-2}, issn={1387-0874}, mr={2892587}}
\bptok{imsref}%
\end{barticle}
\endbibitem

\bibitem{BR}
\begin{barticle}[mr]
\bauthor{\bsnm{Bickel},~\bfnm{Peter~J.}\binits{P.~J.}} \AND
  \bauthor{\bsnm{Ritov},~\bfnm{Ya'acov}\binits{Y.}}
(\byear{2003}).
\btitle{Nonparametric estimators which can be ``plugged-in''}.
\bjournal{Ann. Statist.}
\bvolume{31}
\bpages{1033--1053}.
\bid{doi={10.1214/aos/1059655904}, issn={0090-5364}, mr={2001641}}
\bptok{imsref}%
\end{barticle}
\endbibitem

\bibitem{BM}
\begin{barticle}[mr]
\bauthor{\bsnm{Birg{\'e}},~\bfnm{Lucien}\binits{L.}} \AND
  \bauthor{\bsnm{Massart},~\bfnm{Pascal}\binits{P.}}
(\byear{1995}).
\btitle{Estimation of integral functionals of a density}.
\bjournal{Ann. Statist.}
\bvolume{23}
\bpages{11--29}.
\bid{doi={10.1214/aos/1176324452}, issn={0090-5364}, mr={1331653}}
\bptok{imsref}%
\end{barticle}
\endbibitem

\bibitem{CDG}
\begin{barticle}[mr]
\bauthor{\bsnm{Cl{\'e}ment},~\bfnm{Emmanuelle}\binits{E.}},
  \bauthor{\bsnm{Delattre},~\bfnm{Sylvain}\binits{S.}} \AND
  \bauthor{\bsnm{Gloter},~\bfnm{Arnaud}\binits{A.}}
(\byear{2013}).
\btitle{An infinite dimensional convolution theorem with applications to the
  efficient estimation of the integrated volatility}.
\bjournal{Stochastic Process. Appl.}
\bvolume{123}
\bpages{2500--2521}.
\bid{doi={10.1016/j.spa.2013.04.004}, issn={0304-4149}, mr={3054534}}
\bptok{imsref}%
\end{barticle}
\endbibitem

\bibitem{JP}
\begin{bbook}[mr]
\bauthor{\bsnm{Jacod},~\bfnm{Jean}\binits{J.}} \AND
  \bauthor{\bsnm{Protter},~\bfnm{Philip}\binits{P.}}
(\byear{2012}).
\btitle{Discretization of Processes}.
\bseries{Stochastic Modelling and Applied Probability}
\bvolume{67}.
\bpublisher{Springer}, \blocation{Heidelberg}.
\bid{doi={10.1007/978-3-642-24127-7}, mr={2859096}}
\bptok{imsref}%
\end{bbook}
\endbibitem

\bibitem{JRe}
\begin{bmisc}[auto:STB|2013/06/05|13:45:01]
\bauthor{\bsnm{Jacod},~\bfnm{J.}\binits{J.}} \AND
  \bauthor{\bsnm{Rei{\ss}},~\bfnm{M.}\binits{M.}}
(\byear{2013}).
\bhowpublished{A remark on the rates of convergence for integrated volatility
  estimation in the presence of jumps. Preprint. Available at
  \arxivurl{arXiv:1209.4173}}.
\bptok{imsref}%
\end{bmisc}
\endbibitem

\bibitem{JR-tech}
\begin{bmisc}[auto:STB|2013/06/05|13:45:01]
\bauthor{\bsnm{Jacod},~\bfnm{J.}\binits{J.}} \AND
  \bauthor{\bsnm{Rosenbaum},~\bfnm{M.}\binits{M.}}
(\byear{2012}).
\bhowpublished{Estimation of volatility functionals: The case of a $\sqrt{n}$
  window. Technical report, Laboratoire de Probabilit\'es et Mod\`eles
  Al\'eatoires, Univ. Pierre et Marie Curie. Available at
  \arxivurl{arXiv:1212.1997}}.
\bptok{imsref}%
\end{bmisc}
\endbibitem

\bibitem{JS}
\begin{bbook}[mr]
\bauthor{\bsnm{Jacod},~\bfnm{Jean}\binits{J.}} \AND
  \bauthor{\bsnm{Shiryaev},~\bfnm{Albert~N.}\binits{A.~N.}}
(\byear{2003}).
\btitle{Limit Theorems for Stochastic Processes},
\bedition{2nd} ed.
\bseries{Grundlehren der Mathematischen Wissenschaften [Fundamental Principles
  of Mathematical Sciences]}
\bvolume{288}.
\bpublisher{Springer}, \blocation{Berlin}.
\bid{mr={1943877}}
\bptok{imsref}%
\end{bbook}
\endbibitem

\bibitem{M1}
\begin{barticle}[mr]
\bauthor{\bsnm{Mancini},~\bfnm{Cecilia}\binits{C.}}
(\byear{2009}).
\btitle{Non-parametric threshold estimation for models with stochastic
  diffusion coefficient and jumps}.
\bjournal{Scand. J. Stat.}
\bvolume{36}
\bpages{270--296}.
\bid{doi={10.1111/j.1467-9469.2008.00622.x}, issn={0303-6898}, mr={2528985}}
\bptok{imsref}%
\end{barticle}
\endbibitem

\bibitem{M2}
\begin{barticle}[mr]
\bauthor{\bsnm{Mancini},~\bfnm{Cecilia}\binits{C.}}
(\byear{2011}).
\btitle{The speed of convergence of the threshold estimator of integrated
  variance}.
\bjournal{Stochastic Process. Appl.}
\bvolume{121}
\bpages{845--855}.
\bid{doi={10.1016/j.spa.2010.12.001}, issn={0304-4149}, mr={2770909}}
\bptok{imsref}%
\end{barticle}
\endbibitem

\bibitem{MZ}
\begin{barticle}[mr]
\bauthor{\bsnm{Mykland},~\bfnm{Per~A.}\binits{P.~A.}} \AND
  \bauthor{\bsnm{Zhang},~\bfnm{Lan}\binits{L.}}
(\byear{2009}).
\btitle{Inference for continuous semimartingales observed at high frequency}.
\bjournal{Econometrica}
\bvolume{77}
\bpages{1403--1445}.
\bid{doi={10.3982/ECTA7417}, issn={0012-9682}, mr={2561071}}
\bptok{imsref}%
\end{barticle}
\endbibitem

\bibitem{V10}
\begin{barticle}[mr]
\bauthor{\bsnm{Vetter},~\bfnm{Mathias}\binits{M.}}
(\byear{2010}).
\btitle{Limit theorems for bipower variation of semimartingales}.
\bjournal{Stochastic Process. Appl.}
\bvolume{120}
\bpages{22--38}.
\bid{doi={10.1016/j.spa.2009.10.005}, issn={0304-4149}, mr={2565850}}
\bptok{imsref}%
\end{barticle}
\endbibitem

\end{thebibliography}
\end{document}